\documentclass[11pt]{article}
\usepackage{amsfonts,amsmath,amsxtra}
\usepackage{latexsym}
\usepackage{amssymb}

\def\hybrid{\topmargin 0pt      \oddsidemargin 0pt
        \headheight 0pt \headsep 0pt
        \textwidth 16.5cm
        \textheight 23cm
        \marginparwidth 0.0in
        \parskip 5pt plus 1pt   \jot = 1.5ex}
\catcode`\@=11
\def\marginnote#1{}
\newcount\hour
\newcount\minute
\newtoks\amorpm
\hour=\time\divide\hour by60 \minute=\time{\multiply\hour by60
\global\advance\minute by-\hour}
\edef\standardtime{{\ifnum\hour<12 \global\amorpm={am}%
        \else\global\amorpm={pm}\advance\hour by-12 \fi
        \ifnum\hour=0 \hour=12 \fi
      \number\hour:\ifnum\minute<10 0\fi\number\minute\the\amorpm}}
\edef\militarytime{\number\hour:\ifnum\minute<10 0\fi\number\minute}

\def\draftlabel#1{{\@bsphack\if@filesw {\let\thepage\relax
   \xdef\@gtempa{\write\@auxout{\string
      \newlabel{#1}{{\@currentlabel}{\thepage}}}}}\@gtempa
   \if@nobreak \ifvmode\nobreak\fi\fi\fi\@esphack}
        \gdef\@eqnlabel{#1}}
\def\@eqnlabel{}
\def\@vacuum{}
\def\draftmarginnote#1{\marginpar{\raggedright\scriptsize\tt#1}}

\def\draft{\oddsidemargin -0.1truein
        \def\@oddfoot{\sl preliminary draft \hfil
        \rm\thepage\hfil\sl\today\quad\militarytime}
        \let\@evenfoot\@oddfoot \overfullrule 3pt
        \let\label=\draftlabel
        \let\marginnote=\draftmarginnote
\def\@eqnnum{{\rm (\theequation)}
\rlap{\kern\marginparsep\tt\@eqnlabel}%
\global\let\@eqnlabel\@vacuum}  }

\newfont{\Bbbb}{msbm7 scaled 1\@ptsize00}
\newcommand{\zs}{\raise-1pt\hbox{$\mbox{\Bbbb Z}$}}

scaled 1\@ptsize00

\font\sevenmsa=msam6 
\newfam\msafam
\textfont\msafam=\sevenmsa
\def\hexnumber@#1{\ifnum#1<10 \number#1\else
\ifnum#1=10 A\else\ifnum#1=11 B\else\ifnum#1=12 C\else \ifnum#1=13
D\else\ifnum#1=14 E\else\ifnum#1=15 F\fi\fi\fi\fi\fi\fi\fi}
\def\msa@{\hexnumber@\msafam}
\def\llcorner{\delimiter"4\msa@78\msa@78 }
\def\lrcorner{\delimiter"5\msa@79\msa@79 }
\mathchardef\blacktriangleright="3\msa@49
\mathchardef\blacktriangleleft="3\msa@4A \font\tenmsb=msbm10 scaled
1\@ptsize00
\newfam\msbfam
\textfont\msbfam=\tenmsb \scriptfont\msbfam=\tenmsb

\newdimen\Squaresize \Squaresize=14pt
\newdimen\Thickness \Thickness=0.5pt

\def\Square#1{\hbox{\vrule width \Thickness
   \vbox to \Squaresize{\hrule height \Thickness\vss
      \hbox to \Squaresize{\hss#1\hss}
   \vss\hrule height\Thickness}
\unskip\vrule width \Thickness} \kern-\Thickness}

\def\Vsquare#1{\vbox{\Square{$#1$}}\kern-\Thickness}

\def\numberbysection{\@addtoreset{equation}{section}
        \def\theequation{\thesection.\arabic{equation}}}
\numberbysection

\renewcommand{\theequation}{\thesection.\arabic{equation}}
\def\titlepage{\@restonecolfalse\if@twocolumn\@restonecoltrue\onecolumn
     \else \newpage \fi \thispagestyle{empty}\c@page\z@
        \def\thefootnote{\fnsymbol{footnote}} }

\def\endtitlepage{\if@restonecol\twocolumn \else  \fi
        \def\thefootnote{\arabic{footnote}}
        \setcounter{footnote}{0}}  
\relax

\hybrid
\makeatletter
\newdimen\normalarrayskip            
\newdimen\minarrayskip               
\normalarrayskip\baselineskip \minarrayskip\jot
\newif\ifold             \oldtrue            \def\new{\oldfalse}
\def\arraymode{\ifold\relax\else\displaystyle\fi}
\def\eqnumphantom{\phantom{(\theequation)}} 
\def\@arrayskip{\ifold\baselineskip\z@\lineskip\z@
     \else
     \baselineskip\minarrayskip\lineskip1\baselineskip\fi}

\def\@arrayclassz{\ifcase \@lastchclass \@acolampacol \or
\@ampacol \or \or \or \@addamp \or
   \@acolampacol \or \@firstampfalse \@acol \fi
\edef\@preamble{\@preamble
  \ifcase \@chnum
     \hfil$\relax\arraymode\@sharp$\hfil
     \or $\relax\arraymode\@sharp$\hfil
     \or \hfil$\relax\arraymode\@sharp$\fi}}

\def\@array[#1]#2{\setbox\@arstrutbox=\hbox{\vrule
     height\arraystretch \ht\strutbox
     depth\arraystretch \dp\strutbox
width\z@}\@mkpream{#2}\edef\@preamble{\halign \noexpand\@halignto
\bgroup \tabskip\z@ \@arstrut \@preamble \tabskip\z@ \cr}%
\let\@startpbox\@@startpbox \let\@endpbox\@@endpbox
  \if #1t\vtop \else \if#1b\vbox \else \vcenter \fi\fi
  \bgroup \let\par\relax
  \let\@sharp##\let\protect\relax
  \@arrayskip\@preamble}
\def\eqnarray{\stepcounter{equation}%
              \let\@currentlabel=\theequation
              \global\@eqnswtrue
              \global\@eqcnt\z@
              \tabskip\@centering              
              \let\\=\@eqncr
              $$%
            \halign to \displaywidth  \bgroup
             \eqnumphantom \@eqnsel
      \hskip\@centering                               
    $\displaystyle  \tabskip\z@ {##}$%
    &\global\@eqcnt\@ne \hskip 2\arraycolsep
         $ \displaystyle  \arraymode{##}$\hfil
    &\global\@eqcnt\tw@ \hskip 2\arraycolsep
         $\displaystyle\tabskip\z@{##}$\hfil
         \tabskip\@centering
    &{##}\tabskip\z@\cr}
\makeatother

\newcommand{\ZZ}{{\mathbb{Z}}}

\def\IC{\mathbb{C}}

\def\IR{\mathbb{R}}
\def\IZ{\mathbb{Z}}
\def\CA {\mathcal{A}}

\def\CD {\mathcal{D}}

\def\CH {\mathcal{H}}

\def\CP {\mathcal{P}}

\def\CT {\mathcal{T}}
\def\CU {\mathcal{U}}
\def\CV {\mathcal{V}}

\def\ch{{\cal H}}

\def\ch{{\rm ch}}

\def\nn{\nonumber}

\newtheorem{te}{Theorem}[section]
\newtheorem{de}{Definition}[section]
\newtheorem{prop}{Proposition}[section]           
\newtheorem{cor}{Corollary}[section]
\newtheorem{lem}{Lemma}[section]
\newtheorem{ex}{Example}[section]
\newtheorem{rem}{Remark}[section]

\newcommand\bqa{\begin{eqnarray}}
\newcommand\eqa{\end{eqnarray}}
\def\be{\begin{eqnarray}\new\begin{array}{cc}}
\def\ee{\end{array}\end{eqnarray}}
\def\nn{\nonumber}

\def\beq{\begin{equation}}
\def\eeq{\end{equation}}
\def\bse{\begin{subequations}}                
\def\ese{\end{subequations}}
\def\bp{\begin{pmatrix}}
\def\ep{\end{pmatrix}}

\def\stack#1#2{\raise0.7pt\hbox{$\mathrel{\mathop{#2}\limits^{#1}}$}}
\def\tr{\triangleright}
\def\tl{\triangleleft}
\def\sem{\mathsurround=0pt \raise1pt
\hbox{$\scriptscriptstyle>\!\!$}\:\!\!\tl}
\def\mes{\mathsurround=0pt \tr\!\:\!\raise0.8pt
\hbox{$\scriptscriptstyle\!\!<$}\,}
\def\]{\mathsurround=0pt ]\raise-2pt\hbox{$_\ast$}}

\def\<{\langle}
\def\>{\rangle}

\def\vr{\varrho}

\def\ch{{\cal H}}

\def\CU{{\cal U}}

\def\CH{\mathcal{H}}

\def\we{\raise-1pt\hbox{$\,\stackrel{\wedge}{,}\,$}}
\def\tr{{\rm tr}\,}

\newcounter{pac}[section]

0

0

\title{\bf On $q$-deformed $\mathfrak{gl}_{\ell+1}$-Whittaker function III}
\begin{document}
\author{Anton Gerasimov, Dimitri Lebedev, and Sergey Oblezin}
\date{}

\maketitle

\renewcommand{\abstractname}{}

\begin{abstract}
\noindent {\bf Abstract}. We identify $q$-deformed
$\mathfrak{gl}_{\ell+1}$-Whittaker functions with a specialization
of Macdonald polynomials. This provides a representation of  $q$-deformed
$\mathfrak{gl}_{\ell+1}$-Whittaker functions in terms of Demazure characters of
affine Lie algebra $\widehat{\mathfrak{gl}}_{\ell+1}$.
 We also define a system of dual Hamiltonians for $q$-deformed
$\mathfrak{gl}_{\ell+1}$-Toda chains and  give a new integral
representation for $q$-deformed $\mathfrak{gl}_{\ell+1}$-Whittaker
functions. Finally an expression of  $q$-deformed
$\mathfrak{gl}_{\ell+1}$-Whittaker function  as a matrix element of
a quantum torus algebra is derived.

\end{abstract}

\section*{Introduction}

In \cite{GLO1} an explicit expression for a $q$-deformed
$\mathfrak{gl}_{\ell+1}$-Whittaker function was proposed.
This expression provides a $q$-version of the
Casselman-Shalika-Shintani formula \cite{Sh}, \cite{CS}.
More precisely the $q$-deformed $\mathfrak{gl}_{\ell+1}$-Whittaker
function is given by a character of an infinite-dimensional
$GL(\ell+1,\IC)\times \IC^*$-module. It was remarked in \cite{GLO1} that multiplied
by a simple factor the $q$-deformed $\mathfrak{gl}_{\ell+1}$-Whittaker
functions  have a representation as character of a {\it
  finite-dimensional} $GL(\ell+1,\IC)\times \IC^*$-modules. In this
note  we identify these modules as particular Demazure modules of
affine Lie algebra $\widehat{\mathfrak{gl}}_{\ell+1}$ (see Theorem
\ref{corSan}). This  easily follows from two interpretations of
Macdonald polynomials $P_{\lambda}(x;q,t)$ specialized at $t=0$.
Below we express $q$-deformed $\mathfrak{gl}_{\ell+1}$-Whittaker functions in terms of
$P_{\lambda}(x;q,t=0)$. On the other hand  a relation between
characters of $\widehat{\mathfrak{gl}}_{\ell+1}$ Demazure modules
and $P_{\lambda}(x;q,t=0)$ was established previously by Sanderson
\cite{San1}. Note that the results of \cite{San1} were generalized
to simply-laced semisimple Lie algebras in \cite{I}.
We are going to consider the generalization
of the constructions of  this note to the simply-laced case elsewhere.

The explicit  expression for $q$-deformed
$\mathfrak{gl}_{\ell+1}$-Whittaker function was derived in
\cite{GLO1} by considering a limit $t\to \infty$ of the Macdonald
polynomials $P_{\lambda}(x;q,t)$. In this paper using the same limit we
derive a set of dual Hamiltonian operators of $q$-deformed
$\mathfrak{gl}_{\ell+1}$-Toda chain. The Whittaker function
constructed in \cite{GLO1} is a common  eigenfunction of these dual
Hamiltonian operators as well as standard Hamiltonian operators of
$q$-deformed $\mathfrak{gl}_{\ell+1}$-Toda chain. We also consider a
limit $t\to 0$ of Macdonald polynomials and relate it with
$q$-deformed $\mathfrak{gl}_{\ell+1}$-Whittaker function. However in
this interpretation of Whittaker function the role of standard
Hamiltonian Toda operators  and the dual ones is reversed. This
leads to a new integral representation of $q$-deformed
$\mathfrak{gl}_{\ell+1}$-Whittaker function which is an analog of
Mellin-Banes integral representation for
$\mathfrak{gl}_{\ell+1}$-Whittaker function \cite{KL}. In some sense
this representation of $q$-deformed Whittaker function is dual to
the one considered in \cite{GLO1}.

According to Kostant \cite{Ko},
$\mathfrak{g}$-Whittaker functions naturally arise as
matrix elements of infinite-dimensional representations of
$\CU(\mathfrak{g})$. Using an embedding of $\CU(\mathfrak{g})$
into a tensor product of several copies of Heisenberg algebras
 one obtains a
realization of $\mathfrak{g}$-Whittaker functions as matrix elements
of several copies of Heisenberg algebras. In this paper we construct
analogous representation of $q$-deformed
$\mathfrak{gl}_{\ell+1}$-Whittaker function as a particular matrix
element of several copies of quantum torus algebras. We demonstrate
that this representation is compatible with a $q$-version of Kostant
representation.

Finally note that  we realize a $q$-deformed Whittaker
function multiplied by simple factor as a character of a
finite-dimensional Demazure module of affine Lie algebra. As for
$q$-deformed Whittaker function {\it per se}  we
describe a representation of $q$-deformed
$\mathfrak{gl}_2$-Whittaker function as a character of  a
certain infinite-dimensional representation in the
cohomology of line bundles over a semi-infinite manifold \cite{GLO1}. This
character can be considered as a proper substitute of a
semi-infinite Demazure character of $\widehat{\mathfrak{gl}}_2$
\cite{GLO2}. We are going to discuss this
interpretation ( and its generalization to $\mathfrak{gl}_{\ell+1}$)
in \cite{GLO3}.

The paper is organized as follows. In Section 1 we describe basic
properties of Macdonald polynomials. In particular, using the
self-duality of Macdonald polynomials we define a dual system of
Macdonald operators. In Section 2 we propose two explicit
expressions for $q$-deformed $\mathfrak{gl}_{\ell+1}$-Whittaker
functions as common eigenfunctions of $q$-deformed
$\mathfrak{gl}_{\ell+1}$-Toda chain.  We also construct a system of
dual Hamiltonians for $q$-deformed $\mathfrak{gl}_{\ell+1}$-Toda
chain. In Section 3  the $q$-deformed
$\mathfrak{gl}_{\ell+1}$-Whittaker functions are identified with
Demazure characters for affine Lie algebra
$\widehat{\mathfrak{gl}}_{\ell+1}$. Finally in Section 4 a
representation  of $q$-deformed $\mathfrak{gl}_{\ell+1}$-Whittaker
function as a matrix element of a quantum torus algebra is derived.

{\em Acknowledgments}: The research of AG was  partly supported by
SFI Research Frontier Programme and Marie Curie RTN Forces Universe
from EU. The research of SO is partially supported by  RF President
Grant MK-134.2007.1. The research was also partially supported by
Grant RFBR-08-01-00931-a. The final part of this work was done while
the first and the second authors were enjoying  the hospitality and
support of the Max-Planck-Institute f\"{u}r Mathematik at Bonn.
It is a great pleasure to thank the Institute for excellent working conditions.

\section{Macdonald polynomials}

In this section we recall the standard facts about Macdonald
polynomials. The basic reference is  \cite{Mac} (see also \cite{Ch}
for details and further developments).

Consider symmetric polynomials in  variables
$(x_1,\ldots,x_{\ell+1})$ over the field $\mathbb{Q}(q,t)$ of
rational functions in $q,t$. Given a partition
$\lambda=(\lambda_1,\lambda_2,\ldots,\lambda_{\ell+1})$, that is the
set of  non-negative integers such that
$\lambda_1\geq\lambda_2\geq\cdots\geq\lambda_{\ell+1}$.  Let $\preceq$ be the
partial ordering on the set of partitions;  precisely, given two
partitions $\lambda',\,\lambda$ we write $\lambda'\preceq\lambda$
when $\lambda'_k\leq\lambda_k$ for $k=1,\ldots,\ell+1$.

Let $m_\lambda$ and $\pi_\lambda$ be  polynomial basis of the space
of symmetric polynomials indexed by partitions $\lambda$:
$$
m_\lambda=\sum_{\sigma\in\mathfrak{S}_{\ell+1}}\,
x_{\sigma(1)}^{\lambda_1}x_{\sigma(2)}^{\lambda_2}\cdot \ldots\cdot
x_{\sigma(\ell+1)}^{\lambda_{\ell+1}},
$$
$$
\pi_\lambda\,=\,\pi_{\lambda_1}\pi_{\lambda_2}\cdot\ldots\cdot
\pi_{\lambda_{\ell+1}},\hspace{1.5cm}\pi_n=\sum_{k=1}^{\ell+1}\,x_k^{n},
$$
where $\mathfrak{S}_{\ell+1}$ is the permutation group of $\ell+1$
elements. Define a scalar product $\<\,,\,\>_{q,t}$ on the space of
symmetric functions over $\mathbb{Q}(q,t)$ as follows
$$
\<\pi_\lambda,\,\pi_{\lambda'}\>_{q,t}\,=\,\delta_{\lambda,\lambda'}\cdot
z_\lambda(q,t),
$$
where
$$
z_\lambda(q,t)\,=\,\prod_{n\geq1}\,n^{m_n}m_n!\cdot
\prod_{\lambda_k\neq 0}
\frac{1-q^{\lambda_k}}{1-t^{\lambda_k}},\hspace{1.5cm}
m_n=\bigl|\{k|\,\lambda_k=n\}\bigr|.
$$
In the following we always imply $q<1$.

\begin{de}\label{Macdef} Macdonald polynomials
 $P_{\lambda}(x;q,t)$
are symmetric polynomial  functions over $\mathbb{Q}(q,t)$ such that\\
$$
P_\lambda\,=\, m_\lambda+\sum_{\lambda'\preceq\,\lambda}
u_{\lambda\lambda'}m_{\lambda'},
$$
with $u_{\lambda\lambda'}\in\mathbb{Q}(q,t)$, and
 for $\lambda\neq\lambda'$
$$
\bigl\<P_\lambda,\, P_{\lambda'}\bigr\>_{q,t}\,=\,0.
$$
\end{de}

In the following we slightly extend the notion of Macdonald polynomials
$P_{\lambda}(x;q,t)$ to the case of generalized  partitions
$\lambda=(\lambda_1,\lambda_2,\ldots,\lambda_{\ell+1})$,
$\lambda_1\geq\lambda_2\geq\cdots\geq\lambda_{\ell+1}$, $\lambda_i\in \IZ$ using the
relation
$$
P_{(\lambda_1,\lambda_2,\ldots,\lambda_{\ell+1})}(x;q,t)=
\Big(\prod_{j=1}^{\ell+1}\,z_j^{\lambda_{\ell+1}}\Big)\,
P_{(\lambda_1-\lambda_{\ell+1},
\lambda_2-\lambda_{\ell+1},\ldots,\lambda_{\ell}-\lambda_{\ell+1},0)}(x;q,t)
$$
Although now $P_{\lambda}(x;q,t)$ are not necessary polynomials we
use the term 'Macdonald polynomial' for thus  defined $P_{\lambda}(x;q,t)$.

Macdonald polynomials can be equivalently characterized as common
eigenfunctions of a set of Hamiltonians $H_r$

 \be H_r\,P_{\lambda}(x;q,t)=c_r(q^{\lambda})\,P_{\lambda}(x;q,t), \ee \be
c_r(q^{\lambda})=\chi_r(q^{\lambda}t^{\vr})= \sum_{I_r}\prod_{i\in I_r}\,
q^{\lambda_{i}}\,\, t^{\vr_i}, \ee where the eigenvalues $\chi_r(z)$ are
characters of fundamental representations $\bigwedge^r\IC^{\ell+1}$
of $\mathfrak{gl}_{\ell+1}$, $\vr_i=\ell+1-i$  and we define
$q^{\lambda}t^{\varrho}=(q^{\lambda_1}
t^{\varrho_1},\ldots,q^{\lambda_{\ell+1}}t^{\varrho_{\ell+1}})$.
Here the sum is over ordered subsets
$$
I_r=\{i_1<i_2<\ldots<i_r\}\subset\{1,2,\ldots,\,\ell+1\}.
$$ Explicitly $H_r$ are given by \be\label{MacdonaldHamiltonian}
H_r=\sum_{I_r} \, t^{r(r-1)/2}\, \prod_{i\in I_r,\,j\notin I_r}
\frac{tx_i-x_j}{x_i-x_j}\,\,\,\,\prod_{i\in I_r} T_{x_i},\qquad
r=1,\ldots, \ell+1, \ee and  difference operators $T_{x_i}$ are
defined as

$$
T_{x_i}f(x_1,\ldots,x_i,\ldots,x_{\ell+1})=f(x_1,\ldots,qx_i,\ldots,x_{\ell+1}),
$$
for $i=1,\ldots,\ell+1.$ The simplest operator is given by \be
H_1=\sum_{i=1}^{\ell+1}\,\prod_{i\neq j}
\frac{tx_i-x_j}{x_i-x_j}\,T_{x_i}. \ee 
Let $t<1$ and 
$$
\Delta(x|q,t)=\prod_{i\neq j}\,\prod_{n=0}^{\infty}
\frac{1-x_ix_j^{-1}q^n}{1-tx_ix_j^{-1}q^n}.
$$
Define  another  scalar
product on symmetric functions of $(\ell+1)$-variables $x_1,\ldots
,x_{\ell+1}$ as follows 
\be\label{newprod}
\<f,g\>'_{q,t}=\frac{1}{(\ell+1)!}\, \,\oint_{\Gamma}
\,\prod_{i=1}^{\ell+1}\,\frac{dx_i}{2\pi\imath x_i} \,
f(x^{-1})\,g(x)\,\Delta(x|q,t), 
\ee
 where the integration domain $\Gamma$  is such that each 
$x_i$ goes around $x_i=0$ and is in
the region defined by inequalities 
$t<|x_i/x_j|<t^{-1}$.  
Difference operators $H^{\mathfrak{gl}_{\ell+1}}_r$ are self-adjoint
with respect to $\<\,,\,\>'_{q,t}$:
$$
\<f,H^{\mathfrak{gl}_{\ell+1}}_{r}\,g\>'_{q,t}\,=\,
\<H^{\mathfrak{gl}_{\ell+1}}_{r}\,f,g\>'_{q,t}\,\,.
$$
The following statement was proved in \cite{AOS}.
\begin{prop}
 The following relations hold

1.\be\label{macrecone}
P^{\mathfrak{gl}_{\ell+1}}_{\lambda}(x;q,t)\,=\, \frac{1}{\ell!}\,\,
\frac{ \<P^{\mathfrak{gl}_\ell}_{\lambda},
P^{\mathfrak{gl}_\ell}_{\lambda}\>_{q,t}}{\<P^{\mathfrak{gl}_\ell}_{\lambda},
P^{\mathfrak{gl}_\ell}_{\lambda}\>'_{q,t}}\times\,\, \\
\times \,\,\int_{\Gamma}
 \prod_{i=1}^{\ell}\,\frac{dy_i}{2\pi\imath y_i}\,\,
C_{\ell+1,\ell}(x,y^{-1}|q,t)
P^{\mathfrak{gl}_{\ell}}_{\lambda}(y;q,t) \,\Delta(y|q,t),\ee where
the integration domain $\Gamma$ is as in \eqref{newprod} 
with the additional conditions  $|x_i y_j^{-1}|<1$, $i=1,\ldots
\ell+1$, $j=1,\ldots ,\ell$.

2. \be P^{\mathfrak{gl}_{\ell+1}}_{\lambda+(\ell+1)^k}(x;q,t)
=\Big(\prod_{j=1}^{\ell+1}\,x_j^{k}\Big)\,
P^{\mathfrak{gl}_{\ell+1}}_{\lambda}(x;q,t).\ee Here
$\lambda+(\ell+1)^k=(\lambda_1+k,\ldots ,\lambda_{\ell}+k,k)$ is a
partition obtained from $\lambda$ by a substitution $\lambda_j\to
\lambda_j+k$, $j=1,\ldots,\ell+1$ and
$$
C_{\ell+1,\ell}(x,y|q,t)=\prod_{i=1}^{\ell+1}\prod_{j=1}^{\ell}
\prod_{n=0}^{\infty} \,\frac{1-tx_iy_jq^n}{1-x_iy_jq^n},
$$
$$
\<P^{\mathfrak{gl}_{\ell}}_\lambda,P^{\mathfrak{gl}_{\ell}}_\lambda\>'_{q,t}\,=\,
\prod_{1\leq i<j\leq\ell}\prod_{n=0}^\infty
\frac{1-t^{j-i}q^{\lambda_i-\lambda_j+n}}
{1-t^{j-i+1}q^{\lambda_i-\lambda_j+n}}\cdot
\frac{1-t^{j-i}q^{\lambda_i-\lambda_j+n+1}}
{1-t^{j-i-1}q^{\lambda_i-\lambda_j+n+1}},
$$

$$
\<P^{\mathfrak{gl}_{\ell}}_\lambda,\,
P^{\mathfrak{gl}_{\ell}}_\lambda\>_{q,t}\,=\,
\prod_{i=1}^{\ell}\,\prod_{k=i}^{\ell}\,\,
\prod_{n=1}^{\lambda_k-\lambda_{k+1}}\,
\frac{1-t^{k-i}\,q^{\lambda_i-\lambda_{k+1}+1-n}}
{1-t^{k+1-i}\,q^{\lambda_i-\lambda_{k+1}-n}},
$$
where  $\lambda_{\ell+1}=0$ is assumed in the last formula.
\end{prop}
These relations provide a recursive construction of Macdonald
polynomials corresponding to  arbitrary partitions.

Macdonald polynomials respect a remarkable symmetry (see e.g.
\cite{Ch}). Let us define the normalized Macdonald polynomial
$\Phi_{\lambda}(x;q,t)$ as follows

 \be\label{renromMac}
\Phi_{\lambda}(x;q,t)\,=\,
t^{\sum_{i=1}^{\ell+1}\lambda_i\rho_i}\prod_{n=0}^\infty\prod_{1\leq
i<j\leq\ell} \frac{1-t^2q^{\lambda_i-\lambda_j+n}}{1-tq^{\lambda_i-\lambda_j+n}}\,\,\,
P_{\lambda}(x;q,t), \ee where $\rho_i=\varrho_i-\ell/2=1-i+\ell/2$ for
$i=1,\ldots,\ell+1$.

 In the following we will always imply that
$t=q^{-k}$, $k\in \IZ$ and $q<1$.
 Then for any partitions $\lambda$ and $\mu$ we have:  \be\label{MacdonaldSelfDuality}
\Phi_{\lambda}(q^{\mu-k\rho};q,t)\,=\,\Phi_{\mu}(q^{\lambda-k\rho};q,t). \ee
Define dual Macdonald Hamiltonians by
\be\label{DualMacdonaldHamiltonian}
H^{\vee}_r(q^{\lambda})\,=\,H_r(q^{\lambda}\, t^{\rho}),\qquad
r=1,\ldots,\ell+1. \ee Normalized Macdonald polynomials satisfy the
following eigenvalue problems.

\begin{prop} For any partitions $\lambda$ and $\mu$ the  normalized
Macdonald polynomials satisfy the following system of equations

 \be\left\{\begin{array}{l}
H_r(x)\,\Phi_{\lambda}(x;q,t)\,=\,c_r(q^{\lambda})\,\Phi_{\lambda}(x;q,t),\\
H^\vee_r(q^\lambda)\,\Phi_{\lambda}(x;q,t)\,=\,c^\vee_r(x)\,
\Phi_{\lambda}(x;q,t),\end{array}\right. \ee
 where
\be
 c_r(q^{\lambda})\,=\,\chi_{r}(q^{\lambda}t^{\varrho})\,=\,
\sum_{I_r}\prod_{i\in I_r}\,q^{\lambda_i}\,t^{\varrho_i},\\
c^\vee_r(x)\,=\,\chi_{r}(x\,t^{\ell/2})\,=\,t^{r\ell/2}\,
\sum_{I_r}\prod_{i\in I_r}\,x_i. \ee
\end{prop}
\noindent {\it Proof:} Let $\mu$ be any partition and let
$x=q^{\mu},$ then
$$
H_r({q^{\mu}})\,\Phi_{\lambda}(q^{\mu};q,t)\,=t^{\frac{r\ell}{2}}\,\sum_{I_r}\Big(\prod_{i\in
I_r}\, q^{\lambda_i}\,t^{\rho_i}\Big)\,\Phi_{\lambda}(q^{\mu};q,t).
$$
Let us make a change  variables $\mu\rightarrow\mu-k\rho$.  Then
using  self-duality (\ref{MacdonaldSelfDuality}) of Macdonald
polynomials  one obtains \be
H_r({q^{\mu}t^{\rho}})\,\Phi_{\mu}(q^{\lambda}t^{\rho};q,t)\,=\,
t^{\frac{r\ell}{2}}\,\sum_{I_r}\Big(\prod_{i\in I_r}\,
q^{\lambda_i}\,t^{\rho_i}\Big)\,\Phi_{\mu}(q^{\lambda}t^{\rho};q,t).\ee
Shifting variables $\lambda\rightarrow\lambda+k\rho$ we have \be
H_r({q^{\mu}t^{\rho}})\,\Phi_{\mu}(q^{\lambda};q,t)\,=\,
t^{\frac{r\ell}{2}}\sum_{I_r}\Big(\prod_{i\in I_r}\,
q^{\lambda_i}\,\Big)\,\Phi_{\mu}(q^{\lambda};q,t).\ee
Note that $\Phi_{\lambda}(x;q,t)$ are polynomials in $x$ and thus can be
characterized by its values at $x=q^{\mu}$, $\mu\in \IZ^{\ell+1}$.
Interchanging variables $\lambda\leftrightarrow\mu$  and   denoting
$H^{\vee}_r(q^{\lambda})\,=\,H_r({q^{\lambda}t^{\rho}})$ we obtain
the statement of the proposition $\Box$

\section{ $q$-deformed  $\mathfrak{gl}_{\ell+1}$-Whittaker
function}

In \cite{GLO1} an explicit construction of a $q$-deformed
$\mathfrak{gl}_{\ell+1}$-Whittaker function
$\Psi^{\mathfrak{gl}_{\ell+1}}_{\underline{z}}(\underline{p}_{\ell+1})$
on the lattice
$\underline{p}_{\ell+1}=(p_{\ell+1,1},\ldots,p_{\ell+1,\ell+1})
\in\ZZ^{\ell+1}$ was given. The construction is based on a
particular degeneration of the defining relations for Macdonald
polynomials. In this section using the same degeneration we define
dual Hamiltonians for $q$-deformed $\mathfrak{gl}_{\ell+1}$-Toda
chain. We also consider another degeneration procedure which also
leads to $q$-deformed Toda chain but the role of the Hamiltonians
and the dual Hamiltonians is interchanged. This leads to the second
explicit expression for $q$-deformed Whittaker functions considered
as common eigenfunctions of (dual) Hamiltonians of $q$-deformed Toda
chain.

\subsection{First explicit formula}

The $q$-deformed $\mathfrak{gl}_{\ell+1}$-Whittaker functions are a
common eigenfunction of $q$-deformed $\mathfrak{gl}_{\ell+1}$-Toda
chain Hamiltonians:  \be\label{comm}
\CH_r^{\mathfrak{gl}_{\ell+1}}(\underline{p}_{\ell+1})\,=\,\sum_{I_r}\,\bigl(
\widetilde{X}_{i_1}^{1-\delta_{i_2-i_1,\,1}}\cdot\ldots\cdot
\widetilde{X}_{i_{r-1}}^{1-\delta_{i_r-i_{r-1},\,1}}\cdot
\widetilde{X}_{i_r}^{1-\delta_{i_{r+1}-i_r,\,1}}\bigr)
T_{i_1}\cdot\ldots\cdot T_{i_r},\ee where we assume
$i_{r+1}=\ell+2$.  We use here the following notations
$$
T_if(\underline{p}_{\ell+1})=f(\underline{\widetilde{p}}_{\ell+1})
\hspace{1.5cm}\widetilde{p}_{\ell+1,k}=p_{\ell+1,k}+\delta_{k,i},
$$
and
$$
\widetilde{X}_i=1-q^{p_{\ell+1,i}-p_{\ell+1,i+1}+1},\hspace{0.5cm}
i=1,\ldots,\ell\hspace{1.5cm}\widetilde{X}_{\ell+1}=1.
$$
The   first nontrivial Hamiltonian  is given by:
\be\label{FirstHamiltonian}
 {\CH}_1^{\mathfrak{gl}_{\ell+1}}(\underline{p}_{\ell+1})\,=\,
\sum\limits_{i=1}^{\ell}(1-q^{p_{\ell+1,i}-p_{\ell+1,i+1}+1})T_i\,+\,
T_{\ell+1}. \ee
The corresponding eigenvalue problem can be written in the following
form:
 \be\label{eiglat}
\ch_r^{\mathfrak{gl}_{\ell+1}}(\underline{p}_{\ell+1})
\Psi^{\mathfrak{gl}_{\ell+1}}_{z_1,\ldots ,z_{\ell+1}}
(\underline{p}_{\ell+1})\,=\,(\sum_{ I_r}\prod\limits_{i\in I_r} z_i
)\,\,\Psi^{\mathfrak{gl}_{\ell+1}}_{z_1,\ldots ,z_{\ell+1}}
(\underline{p}_{\ell+1}).\ee The main result of \cite{GLO1} can be
formulated as follows. Denote by
$\CP^{(\ell+1)}\subset\ZZ^{\ell(\ell+1)/2}$ a  subset of  parameters
$p_{k,i}$, $k=1,\ldots,\ell$, $i=1,\ldots,k$ satisfying the
Gelfand-Zetlin conditions $p_{k+1,i}\geq p_{k,i}\geq p_{k+1,i+1}$.
Let $\CP_{\ell+1,\ell}\subset\CP^{(\ell+1)}$ be a set of
$\underline{p}_\ell=(p_{\ell,1},\ldots,p_{\ell,\ell})$ satisfying
the conditions $p_{\ell+1,i}\geq p_{\ell,i}\geq p_{\ell+1,i+1}$.

\begin{te}
The common solution of the eigenvalue problem (\ref{eiglat}) can be
written in the following form.  For $\underline{p}_{\ell+1}$ being
in the dominant  domain $p_{\ell+1,1}\geq\ldots\geq
p_{\ell+1,\ell+1}$ the solution is given by \be\label{main}
 \Psi^{\mathfrak{gl}_{\ell+1}}_{z_1,\ldots,z_{\ell+1}}
(\underline{p}_{\ell+1})\,=\, \sum_{p_{k,i}\in\CP^{(\ell+1)}}\,\,
\prod_{k=1}^{\ell+1} z_k^{\sum_i  p_{k,i}-\sum_i p_{k-1,i}}\,\,
\\ \times\frac{\prod\limits_{k=2}^{\ell}\prod\limits_{i=1}^{k-1}
(p_{k,i}-p_{k,i+1})_q!}
{\prod\limits_{k=1}^{\ell}\prod\limits_{i=1}^k
(p_{k+1,i}-p_{k,i})_q!\,\, (p_{k,i}-p_{k+1,i+1})_q!},\ee where  we
use the notation $(n)_q!=(1-q)...(1-q^n)$.  When
$\underline{p}_{\ell+1}$ is outside the dominant domain we set
$$
\Psi^{\mathfrak{gl}_{\ell+1}}_{z_1,\ldots,z_{\ell+1}}
(p_{\ell+1,1},\ldots,p_{\ell+1,\ell+1})\,=\,0.
$$
\end{te}
Formula (\ref{main}) can be written in the recursive form.
\begin{cor}
 The following recursive relation holds
\be\label{qtodarec}
 \Psi^{\mathfrak{gl}_{\ell+1}}_{z_1,\ldots,z_{\ell+1}}
(\underline{p}_{\ell+1})\,=\,
\sum_{\underline{p}_\ell\in\CP_{\ell+1,\ell}}\,\,
\Delta(\underline{p}_{\ell}) \,\,z_{\ell+1}^{\sum_i
p_{\ell+1,i}-\sum_i p_{\ell,i}}\,\,
Q_{\ell+1,\ell}(\underline{p}_{\ell+1},\underline{p}_{\ell}|q)
\Psi^{\mathfrak{gl}_{\ell}}_{z_1,\ldots,z_{\ell}}(\underline{p}_{\ell}),
\nonumber \ee where \be
Q_{\ell+1,\ell}(\underline{p}_{\ell+1},\underline{p}_{\ell}|q)\,=\,
\frac{1}{\prod\limits_{i=1}^{\ell} (p_{\ell+1,i}-p_{\ell,i})_q!\,\,
(p_{\ell,i}-p_{\ell+1,i+1})_q!},\\
\Delta(\underline{p}_{\ell})=
\prod_{i=1}^{\ell-1}(p_{\ell,i}-p_{\ell,i+1})_q!\,\,\,. \ee
\end{cor}

\begin{lem} \label{funchar}  The $q$-deformed $\mathfrak{gl}_{\ell+1}$-Whittaker
function at  $p_{\ell+1,i}=k+1$ for $i\leq r$, $p_{\ell+1,i}=k$, for
$i> r$ is proportional  to the  character $\chi_r(z)$ of the
fundamental representation  $\Lambda^{r}\IC$ of
 $\mathfrak{gl}_{\ell+1}$
$$
 \Psi^{\mathfrak{gl}_{\ell+1}}_{z_1,\ldots,z_{\ell+1}}
(k+1,\ldots k+1,k,\ldots k)= \Big(\prod_{i=1}^{\ell+1}z_i^k\Big)\,
\chi_r(z)=
\Big(\prod_{i=1}^{\ell+1}z_i^k\Big)\,\,\sum_{I_r}\prod_{i\in
I_r}\,z_i\, .
$$
\end{lem}
\nn {\it Proof}: Directly follows from the general expression
\eqref{main} $\Box$

\begin{ex} Let $\mathfrak{g}=\mathfrak{gl}_{2}$,
$p_{2,1}:=p_1\in\ZZ$, $ p_{2,2}:=p_2\in\ZZ$ and $p_{1,1}:=p\in \ZZ$.
The function
$$
\Psi_{z_1,z_2}^{{\mathfrak gl}_2}(p_{1},p_{2}) =\sum_{p_{2}\leq
p\leq p_{1}}\frac{ z_1^{p} z_2^{p_{1}+p_{2}-p}}
{(p_1-p)_q!(p-p_2)_q!},\qquad  p_{1}\geq p_{2}\,,
$$
$$
\Psi_{z_1,z_2}^{{\mathfrak gl}_2}(p_{1},p_{2})=0, \qquad
p_{1}<p_{2}\,,
$$
is a  common eigenfunction of  mutually commuting Hamiltonians
$$
{\cal H}_1^{{\mathfrak gl}_2}\,=\,(1-q^{p_1-p_2+1})T_1+T_2,\qquad
{\cal H}_2^{{\mathfrak gl}_2}=T_1T_2.
$$
\end{ex}

\subsection{Dual Hamiltonians for $\mathfrak{gl}_{\ell+1}$-Toda chain}

The Hamiltonian operators of $q$-deformed
$\mathfrak{gl}_{\ell+1}$-Toda chain can be obtained by a
degeneration of Macdonald operators discussed in the previous
section (see e.g. \cite{GLO1}). Similarly the  degeneration of dual
Macdonald operators leads to a set of dual Hamiltonians of
$q$-deformed $\mathfrak{gl}_{\ell+1}$-Toda chain.
\begin{prop}1. Let $t=q^{-k}$, $q<1$. Define the limit
$k\to \infty$ of the Macdonald (dual) operators
 \be\label{commx}\CH_r(x)\,=\,\lim_{k\to
\infty} D(x)\,H_r(xq^{-k\rho})\,D(x)^{-1}\,=\\=\,
\sum_{I_r}\,\bigl(X_{i_1}^{1-\delta_{i_1,\,1}}\cdot
X_{i_2}^{1-\delta_{i_2-i_1,\,1}}\cdot\ldots\cdot
X_{i_r}^{1-\delta_{i_r-i_{r-1},\,1}} \bigr)
T_{x_{i_1}}\cdot\ldots\cdot T_{x_{i_r}},\ee
 \be
\CH_r^{\vee}(q^{\lambda})\,=\,\lim_{k\to\infty}\,
q^{kr(2\ell+1-r)/2}\,G(q^\lambda)\,
H^\vee_r(q^{\lambda+k\varrho})\,G(q^\lambda)^{-1}\,=\\=\,\
\,q^{r(r-1)/2}\,\sum_{I_r}\prod_{i\in I_r,\,j\notin I_r}\,
\frac{q^{\lambda_j}}{q^{\lambda_j}-q^{\lambda_i}} \,\,\,\prod_{i\in I_r}
T_{\lambda_i},\ee
here $r=1,\ldots,\ell+1$ and we set $X_i=1-x_{i-1}^{-1}x_i$,
$X_1=1$, $T_{\lambda_i}\lambda_j=\lambda_jT_{\lambda_i}+\delta_{ij}$
 and we assume
$$
D(x)=\prod_{i=1}^{\ell+1}x_i^{-k\varrho_i},
$$
\be G(q^\lambda)
\,=\,(-1)^{\ell\sum_{i=1}^{\ell+1}\lambda_i}\,\,
q^{-\ell\sum_{i=1}^{\ell+1}\lambda_i/2}\,
\prod_{i<j}\,q^{(\lambda_i-\lambda_j)^2/2}.\ee
2. Let \be \Psi_{\lambda}(x)=
\lim_{k\to\infty}\,\,G(q^\lambda)D(x)\,\,
\Phi_{\lambda+k\vr}(xq^{-k\rho};q,t), \ee
 then the following relations
hold \be\CH_r(x)\,\Psi_{\lambda}(x)\,=\,
\chi_r(q^\lambda)\, \Psi_{\lambda}(x),\\
\CH^\vee_r(q^\lambda)\,\Psi_{\lambda}(x)\,=\,
\Big(q^{r(r-1)/2}\prod_{i=1}^r\,x_{i}\Big)\,\Psi_{\lambda}(x),\ee for
$r=1,\ldots,\ell+1$.
\end{prop}

\noindent {\it Proof}: Direct calculations $\Box$

 Observe that the following relation between (\ref{commx}) and
(\ref{comm}) holds
$$
\CH^{\mathfrak{gl}_{\ell+1}}(\underline{p}_{\ell+1})\,=\,\CH_r(x),\hspace{1.5cm}
x_i=q^{p_{\ell+1,\,\ell+2-i}+\vr_{\ell+2-i}},\qquad i=1,\ldots,\ell+1,
$$
for $r=1,\ldots,\ell+1$.

The limit $t=q^{-k}\to \infty$, $k\to \infty$ was used in \cite{GLO1} to
obtain Hamiltonians of $q$-deformed $\mathfrak{gl}_{\ell+1}$-Toda
chain. There is a ``dual'' limit $t=q^{-k}$, $k\to -\infty$ which
also leads to  (dual) Hamiltonians of $q$-deformed
$\mathfrak{gl}_{\ell+1}$-Toda chain but the Hamiltonians and the
dual Hamiltonians are interchanged.

\begin{prop}
1. Let $t=q^{-k}$, $q<1$. Define the limit $k\to -\infty$ of the
Macdonald (dual) operators and their common eigenfunction  as
follows

\be\widehat{\CH}_r(x)\,=\,\lim_{k\to-\infty}\,q^{kr(r-1)/2}\,
H_r(x)\,=\, \sum_{I_r}\prod_{i\in I_r,\,j\notin I_r}\,
\frac{x_j}{x_j-x_i} \,\,\,\prod_{i\in I_r} T_{x_i},\ee

\be\widehat{\CH}_r^{\,\vee}(q^\lambda)\,=\,\lim_{k\to-\infty}\,\,
q^{kr\ell/2}\,\widehat{D}(q^\lambda)\,H^\vee_r(q^{\lambda})\,
\widehat{D}(q^\lambda)^{-1}\,=\\=\,\sum_{I_r}\,\bigl(
\widehat{X}_{i_1}^{1-\delta_{i_1,\,1}}\cdot
\widehat{X}_{i_2}^{1-\delta_{i_2-i_1,\,1}}\cdot\ldots\cdot
\widehat{X}_{i_r}^{1-\delta_{i_r-i_{r-1},\,1}} \bigr)
T_{\lambda_{i_1}}\cdot\ldots\cdot T_{\lambda_{i_r}}, \ee
$\widehat{X}_i\,=\,1-q^{\lambda_i-\lambda_{i+1}}\,$ and
$\widehat{X}_1=1$. We assume here
$\widehat{D}(q^\lambda)=\prod_{i=1}^{\ell+1}q^{k\lambda_i\rho_i}$.

2.  Define
\be\label{Mact0}\widehat{\Psi}_\lambda(x)\,=\,\lim_{k\to-\infty}
\widehat{D}(q^\lambda)\,\Phi_\lambda(x;q,t).\ee Then the following
equations hold
\be\label{modpsi}\widehat{\CH}_r(x)\,\widehat{\Psi}_\lambda(x)\,=\,
q^{\lambda_{\ell+2-r}+\ldots+\lambda_{\ell+1}}\,\widehat{\Psi}_\lambda(x),\\
\widehat{\CH}^{\,\vee}_r(q^\lambda)\,\widehat{\Psi}_\lambda(x)\,=\,
\chi_r(x)\,\widehat{\Psi}_\lambda(x).\ee
\end{prop}

\emph{Proof. } 1. The formula for $\widehat{\CH}_r$ follows straightforwardly.
2. For $t=q^{-k}$ we obtain
$$
\widehat{D}(q^\lambda)\,H^\vee_r(q^{\lambda})\,\widehat{D}(q^\lambda)^{-1}\,=\,
t^{r\ell/2}\sum_{I_r}\,\,\prod_{i\in I_r,\,j\notin I_r}\prod_{i<j}
\frac{t^{j+1-i}q^{\lambda_i}-q^{\lambda_j}}{t^{j-i}q^{\lambda_i}-q^{\lambda_j}}
\prod_{i>j}
\frac{q^{\lambda_i}-t^{i-1-j}q^{\lambda_j}}{q^{\lambda_i}-t^{i-j}q^{\lambda_j}}\,\,
\prod_{i\in I_r}T_{\lambda_i},
$$
due to the following identity
$$
\frac{r(r-1)}{2}\,+\,\sum_{i\in I_r}\bigl(\rho_i+b_{i,\,I_r}\bigr)
\,=\,\frac{r\ell}{2},
$$
where $b_{i,I_r}=\bigl|\bigl\{j\notin I_r\,|\,j<i\bigr\}\bigr|$.

 Thus under the limit $t\to0$ one gets the following.
$$
\frac{t^{j+1-i}q^{\lambda_i}-q^{\lambda_j}}{t^{j-i}q^{\lambda_i}-q^{\lambda_j}}
\longrightarrow1,\,\,i<j, \hspace{1cm}
\frac{q^{\lambda_{i+1}}-q^{\lambda_i}}{q^{\lambda_{i+1}}-tq^{\lambda_i}}
\longrightarrow1-q^{\lambda_i-\lambda_{i+1}}, \hspace{1cm}
\frac{q^{\lambda_i}-t^{i-1-j}q^{\lambda_j}}{q^{\lambda_i}-t^{i-j}q^{\lambda_j}}
\longrightarrow1,\,\,i>j+1,
$$
$\Box$

\begin{rem}
Let $\lambda$ be a partition, then

1.\be\widehat{\CH}_r(x_1,\ldots,x_{\ell+1})\,=\,q^{-\frac{r(r-1)}{2}}
\CH^{\vee}_r(x_1,\ldots,{x_{\ell+1}}),\\
\widehat{\CH}^{\,\vee}_r(q^{\lambda_1},\ldots,q^{\lambda_{\ell+1}})\,=\,
\Delta(q^\lambda)\,
\CH_r(q^{\lambda_{\ell+1}+\vr_{\ell+1}},\ldots,q^{\lambda_1+\vr_1})\,
\Delta(q^\lambda)^{-1}, \ee for $r=1,\ldots,\ell+1$ and
\be\Delta(q^\lambda)\,=\,\prod_{i=1}^{\ell}\,(\lambda_i-\lambda_{i+1})_q!\ee
2. The specialization of Macdonald polynomial at $t=0$\be
\widehat{\Psi}_{\lambda}(x)=P_{\lambda}(x;q,t=0),\ee satisfies equations (\ref{modpsi}).
\end{rem}

\nn {\it Proof}:   Proof of (1) is straightforward and the statement of
(2) easily follows from \eqref{Mact0} and \eqref{renromMac} $\Box$

\subsection{Second explicit formula }

Now we construct an integral representation for $q$-deformed
Whittaker functions by taking $t\to 0$ limit of the recursive
construction of Macdonald polynomials.

In the limit $t\to 0$  the Macdonald scalar product on symmetric
functions of $(\ell+1)$-variables $x_1,\ldots ,x_{\ell+1}$ is
reduced to
\be\label{newprod0}
\<f,g\>'_{q,t=0}=\frac{1}{(\ell+1)!}\oint_{\Gamma_0}
\,\prod_{i=1}^{\ell+1}\,\frac{dx_i}{2\pi\imath x_i} \,
f(x^{-1})\,g(x)\,\Delta(x|q,t=0),
\ee
where
$$
\Delta(x|q,0)=\prod_{i\neq j}\,\prod_{n=0}^{\infty}
(1-x_ix_j^{-1}q^n).
$$
and  the integration domain $\Gamma_0$  is such that each 
$x_i$ goes over a small circle around $x_i=0$.

The limit $t\to 0$ of the recursive kernel $C_{\ell+1,\ell}$ is
given by
$$
C_{\ell+1,\ell}(x,y|q,t=0)=\prod_{i=1}^{\ell+1}\prod_{j=1}^{\ell}
\prod_{n=0}^{\infty} \,\frac{1}{1-x_iy_jq^n}.
$$

\begin{prop}\label{reclimt0}
 \begin{enumerate}
\item Given a partition
 $\lambda=(\lambda_1,\ldots,\lambda_{\ell})$ the following
recursive relation holds \be\label{macrec}
P^{\mathfrak{gl}_{\ell+1}}_{\lambda}(x;q,t=0)\,=\,
\frac{A_{\ell}}{\ell!}\,\,\,\int_{\Gamma_0}\,\,\,\prod_{i=1}^{\ell}\frac{dy_i}{2\pi\imath
y_i}\,\,\, C_{\ell+1,\ell}(x,y^{-1}|q,t=0)\,\times
\\ \times P^{\mathfrak{gl}_{\ell}}_{\lambda}(y;q,t=0)
\,\Delta(y|q,t=0), \ee where
$$
A_{\ell}=\lim_{t\to 0}\,\,\frac{ \<P^{\mathfrak{gl}_\ell}_{\lambda},
P^{\mathfrak{gl}_\ell}_{\lambda}\>_{q,t}}{\<P^{\mathfrak{gl}_\ell}_{\lambda},
P^{\mathfrak{gl}_\ell}_{\lambda}\>'_{q,t}}\,=\,
\prod_{m=1}^\infty(1-q^m)^{\ell-1}\cdot (\lambda_{\ell})_q!\,,
$$and  the contour of integration $\Gamma_0$ is as in \eqref{newprod0}
with additionally conditions $x_iy_j^{-1}<1$.
\item Given a partition
 $\lambda=(\lambda_1,\ldots,\lambda_{\ell})$
 \be P^{\mathfrak{gl}_{\ell+1}}_{\lambda+(\ell+1)^k}(x;q,t=0)
=\Big(\prod_{j=1}^{\ell+1}\,x_j^{k}\Big)\,
P^{\mathfrak{gl}_{\ell+1}}_{\lambda}(x;q,t=0), \ee
where $\lambda+(\ell+1)^k=(\lambda_1+k,\ldots ,\lambda_{\ell}+k,k)$.
\end{enumerate}
\end{prop}

\nn {\it Proof}: We have
$$
\<P^{\mathfrak{gl}_{\ell}}_\lambda,\,
P^{\mathfrak{gl}_{\ell}}_\lambda\>'_{q,t=0}\,=\,
\prod_{i=1}^{\ell-1}\prod_{m=1}^{\infty}\,
\frac{1}{(1-q^{\lambda_i-\lambda_{i+1}+m})},
$$
$$
\<P^{\mathfrak{gl}_{\ell}}_\lambda,\,
P^{\mathfrak{gl}_{\ell}}_\lambda\>_{q,t=0}\,=\,\prod_{i=1}^{\ell-1}
(\lambda_{i}-\lambda_{i+1})_q!\times (\lambda_{\ell})_q!
$$
where   $\prod_{m=1}^{0}(1-q^{m})=1$ is assumed. Thus we obtain the
recursive relation \eqref{macrec} $\Box$

These relations provide a recursive construction of a  $q$-deformed
Whittaker function corresponding to  an arbitrary partition. Note
that the property of Macdonald polynomial being symmetric function
of variables $z_1,\ldots , z_{\ell+1}$  remains true in the limit
$t\to 0$.

\begin{prop}  Let $z_i:=x_{\ell+1,i}$ for $i=1,\ldots,\ell+1$.
Define the function ${}^{MB}\Psi^{\mathfrak{gl}_{\ell+1}}
_{z_1,\ldots, z_{\ell+1}} (\underline{p}_{\ell+1})$ given for the
dominant domain $p_{\ell+1,1}\geq \ldots \geq p_{\ell,1,\ell+1}$ by
an integral expression
 \be {}^{MB}\Psi^{\mathfrak{gl}_{\ell+1}} _{z_1,\ldots, z_{\ell+1}}
(\underline{p}_{\ell+1})\,=\,
(q,q)_{\infty}^{\ell(\ell-1)/2}\,\,\,\int_{\cal{S}}\,\prod_{
n=1;j\leq n}^{\ell}
\frac{dx_{nj}}{2\pi\imath x_{nj}}\\
\times\prod\limits_{n=1}^{\ell+1}\,\prod\limits_{j=1}^{n}
\Big(\frac{x_{n,j}}{x_{n-1,j}}\Big)^{p_{\ell+1,\,n}}
\prod_{n=1}^{\ell}
\frac{\prod\limits_{k=1}^n\prod\limits_{m=1}^{n+1} \Gamma_q(x_{n
k}^{-1}x_{n+1,m}^{})} {n!\,\prod\limits_{s\neq p}
\Gamma_q(x_{ns}x_{np}^{-1})}, \ee where the contour $\cal{S}$ is
obtained by induction from the contours $\Gamma_0$ defined in the
Proposition \ref{reclimt0}  and outside of the dominant domain   by
$$
{}^{MB}\Psi^{\mathfrak{gl}_{\ell+1}}_{z_1,\ldots,z_{\ell+1}}
(p_{\ell+1,1},\ldots,p_{\ell+1,\ell+1})\,=\,0.
$$
Then the function ${}^{MB}\Psi^{\mathfrak{gl}_{\ell+1}}
_{z_1,\ldots, z_{\ell+1}} (\underline{p}_{\ell+1})$ possess the
following properties\\
1. It is $\mathfrak{S}_{\ell+1}$-symmetric:
$$
{}^{MB}\Psi^{\mathfrak{gl}_{\ell+1}} _{z_{\sigma(1)},\ldots,
  z_{\sigma(\ell+1)}}
(\underline{p}_{\ell+1})= {}^{MB}\Psi^{\mathfrak{gl}_{\ell+1}}
_{z_1,\ldots, z_{\ell+1}} (\underline{p}_{\ell+1}), \qquad \sigma\in
\mathfrak{S}_{\ell+1},
$$
2. It is a common   eigenfunction of (dual) Hamiltonians $\CH_r$,
$\CH_r^{\vee}$:
 \be\CH_r^{\mathfrak{gl}_{\ell+1}}(\underline{p}_{\ell+1})\,{}^{MB}
\Psi_{z}(\underline{p}_{\ell+1})\,=\,
\chi_r(z)\, {}^{MB}\Psi_{z}(\underline{p}_{\ell+1}),\\
q^{-r(r-1)/2}\CH^\vee_r(z){}^{MB}\Psi_{z}(\underline{p}_{\ell+1})\,\,=
\Big(\prod_{i=1}^r\, q^{p_{\ell+1,i}}\Big)\,
{}^{MB}\Psi_{z}(\underline{p}_{\ell+1})\,,\ee
for $r=1,\ldots,\ell+1$.

\end{prop}

This integral representation is a $q$-version of Mellin-Barnes
integral representation for $\mathfrak{gl}_{\ell+1}$-Whittaker
functions introduced in \cite{KL}. Let us compare
${}^{MB}\Psi^{\mathfrak{gl}_{\ell+1}} _{{z}}$ with the function
$\Psi^{\mathfrak{gl}_{\ell+1}} _{{z}}$ given by \eqref{main}.

\begin{prop} $q$-deformed $\mathfrak{gl}_{\ell+1}$-Whittaker function
given by \eqref{main}  is a $\mathfrak{S}_{\ell+1}$-symmetric
function
$$
\Psi^{\mathfrak{gl}_{\ell+1}} _{z_{\sigma(1)},\ldots,
  z_{\sigma(\ell+1)}}
(\underline{p}_{\ell+1})= \Psi^{\mathfrak{gl}_{\ell+1}}
_{z_1,\ldots, z_{\ell+1}} (\underline{p}_{\ell+1}), \qquad \sigma\in
\mathfrak{S}_{\ell+1}.
$$
\end{prop}

\nn {\it Proof}: We prove this statement by the induction. Given a
$\mathfrak{gl}_{\ell}$-Whittaker function which is symmetric
$$
\Psi^{\mathfrak{gl}_{\ell}}_{z_{\sigma(1)},\ldots z_{\sigma(\ell)}}
(\underline{p}_{\ell}) =\Psi^{\mathfrak{gl}_{\ell}}_{z_1,\ldots,
  z_{\ell}}(\underline{p}_{\ell}),\qquad \sigma\in \mathfrak{S}_{\ell}.
$$
The function $\Psi^{\mathfrak{gl}_{\ell+1}}$ then given by
$$
\Psi^{\mathfrak{gl}_{\ell}}_{z_1,\ldots,
  z_{\ell},z_{\ell+1}}(\underline{p}_{\ell+1})=
\sum_{\underline{p}_{\ell}\in \CP_{\ell+1,\ell}}
C_{\ell+1,\ell}(q)\,\, z_{\ell+1}^{\sum_{j=1}^{\ell+1}
p_{\ell+1,j}-\sum_{j=1}^{\ell}p_{\ell,j}}
\,\,\Psi^{\mathfrak{gl}_{\ell}}_{z_1,\ldots,
z_{\ell}}(\underline{p}_{\ell}).
$$
The space of solutions of $q$-deformed $\mathfrak{gl}_{\ell+1}$-Toda
chain invariant with respect to $\mathfrak{S}_{\ell}\subset
\mathfrak{S}_{\ell+1}$ is $(\ell+1)$-dimensional. Thus to verify
that \eqref{main} is $\mathfrak{S}_{\ell+1}$-invariant one should
check that it is invariant at $\ell+1$ particular values of
$\underline{p}_{\ell+1}$. Let us take $\underline{p}_{\ell+1}$
corresponding to fundamental representations. By Lemma \ref{funchar}
the corresponding $q$-Whittaker functions are given by characters of
$\mathfrak{gl}_{\ell+1}$-fundamental representations and thus
explicitly $\mathfrak{S}_{\ell+1}$-invariant $\Box$

The function ${}^{MB}\Psi^{\mathfrak{gl}_{\ell+1}}_{z_1,\ldots,z_{\ell+1}}
(\underline{p}_{\ell+1})$ satisfies the full set of equations (i.e.
including dual Hamiltonians) and the function
$\Psi^{\mathfrak{gl}_{\ell+1}}_{z_1,\ldots,z_{\ell+1}}
(\underline{p}_{\ell+1})$  satisfies the original
$q$-deformed Toda equations. Thus one has
$$
\Psi^{\mathfrak{gl}_{\ell+1}} _{z_1,\ldots, z_{\ell+1}}
(\underline{p}_{\ell+1})=C(z_1,\ldots, z_{\ell+1})\,\,\,
{}^{MB}\Psi^{\mathfrak{gl}_{\ell+1}} _{z_1,\ldots, z_{\ell+1}}
(\underline{p}_{\ell+1}).
$$
\begin{prop}
$$
\Psi^{\mathfrak{gl}_{\ell+1}} _{z_1,\ldots, z_{\ell+1}}
(\underline{p}_{\ell+1})\,=\,{}^{MB}\Psi^{\mathfrak{gl}_{\ell+1}}
_{z_1,\ldots, z_{\ell+1}} (\underline{p}_{\ell+1}).
$$
\end{prop}
\noindent{\it Proof:} Denote
$$\widetilde{\Psi}^{\mathfrak{gl}_{\ell+1}} _{z_1,\ldots,
z_{\ell+1}}(\underline{p}_{\ell+1})\,=\,
\Delta(\underline{p}_{\ell+1})\Psi^{\mathfrak{gl}_{\ell+1}}
_{z_1,\ldots, z_{\ell+1}} (\underline{p}_{\ell+1}),$$
$$
\Delta(\underline{p}_{\ell+1}){}^{MB}\Psi^{\mathfrak{gl}_{\ell+1}}
_{z_1,\ldots,
z_{\ell+1}}(\underline{p}_{\ell+1})=P_{\underline{p}_{\ell+1}}(z;q,t=0).$$
Then $ \widetilde{\Psi}^{\mathfrak{gl}_{\ell+1}} _{z_1,\ldots,
z_{\ell+1}} (\underline{p}_{\ell+1})|_{p_{\ell+1,i=0}}\,=\,1 $ and $
P_{(0,0,\cdots,0)}(z;q,t=0)=1 $ by definition of Macdonald
polynomials.  Thus $C(z_1,\ldots, z_{\ell+1})=1$ $\Box$
\begin{rem}\label{RemMac}
The  normalized $q$-deformed $\mathfrak{gl}_{\ell+1}$-Whittaker
function coincides with a $t=0$ specialization of Macdonald
polynomial
\be\widetilde{\Psi}^{\mathfrak{gl}_{\ell+1}}_{z}(\underline{p}_{\ell+1})=
P_{\lambda}(z;q,t=0),\qquad \lambda=(p_{\ell+1,1},\ldots ,p_{\ell+1,\ell+1}).\ee
\end{rem}

\section{$q$-Whitaker functions as characters of affine  Demazure modules}

In this Section we identify the normalized $q$-deformed
$\mathfrak{gl}_{\ell+1}$-Whittaker function $
\widetilde{\Psi}_{z}^{\mathfrak{gl}_{\ell+1}}(\underline{p}_{\ell+1})$
with characters of
 affine Lie algebra $\widehat{\mathfrak{gl}}_{\ell+1}$
 Demazure modules. This straightforwardly follows
from a characterization of the normalized $q$-deformed
$\widehat{\mathfrak{gl}}_{\ell+1}$-Whittaker function as a
specialization of Macdonald polynomials $P_{\lambda}(z;q,t)$ at
$t=0$ (see Remark \ref{RemMac}) and a relation of
$P_{\lambda}(z;q,t=0)$ with  characters of  Demazure modules of
affine Lie algebra $\widehat{\mathfrak{gl}}_{\ell+1}$ established in
\cite{San1}.

To state precisely  the relation between Whittaker functions and Demazure modules
let us start recalling  the notion of a Demazure module \cite{De} (see
 \cite{Ku}, \cite{M} for a general  case of Kac-Moody algebras).
Let $\mathfrak{g}$ be a  Kac-Moody algebra with  Cartan matrix
$\|a_{ij}\|$, $\mathfrak{h}\subset \mathfrak{b}\subset \mathfrak{g}$
 be a Cartan and Borel subalgebras. Let $R\subset \mathfrak{h}^*$ be a
 corresponding root system, $R_+\subset R$ be a subset of
 positive roots corresponding to the Borel subalgebra $\mathfrak{b}$,
$\alpha_1,\ldots ,\alpha_{r}\in R_+$ be a set of simple roots. Denote
$(\lambda,\mu)$ the scalar  product on $\mathfrak{h}^*$ induced by the Killing
form on $\mathfrak{g}$. Given a root $\alpha$ let $\alpha^{\vee}=2\alpha/(\alpha,\alpha)$ be the
corresponding coroot where we identify $\mathfrak{h}\equiv
\mathfrak{h}^*$ using quadratic form $(\,,\,)$. The weight lattice $P$ is
given by $P=\{\lambda\in  \mathfrak{h}^*\,:\,\,(\lambda,\alpha^{\vee}),\in \IZ\,\,\,\,
\alpha\in R\}$. The weight lattice is generated by fundamental weights
$\omega_1,\ldots ,\omega_r$ defined by the conditions
$(\omega_i,\alpha^{\vee}_j)=\delta_{ij}$. The set of dominant weights
is given by $P^+=\{\lambda\in P\,:\,(\lambda,\alpha^{\vee})\geq 0, \,\,\,\,\alpha\in R\}$.
The Weyl group  $W$ is defined as a group of reflections
$s_{\alpha}:\mathfrak{h}^*\to \mathfrak{h}^*$, $\alpha\in R$
$$
s_{\alpha}\,:\,\lambda\,\,\longrightarrow  \lambda-(\lambda,\alpha^{\vee})\alpha,
$$
and is generated by reflections $s_{i}$ corresponding to
simple roots $\alpha_i$. An expression of a Weyl group element $w$ as a
product $w=s_{i_1}\cdots s_{i_l}$ which has minimal length is called
reduced decomposition for $w$ and its length $l(w)=l$ is called a
length of $w$. Let $T$ be a Cartan  torus ${\rm Lie}(T)=\mathfrak{h}$.
The group of characters $X=X(T)$ of $T$ is isomorphic to the  weight
lattice  $P$ of $\mathfrak{g}$. Its group algebra $\IZ[T]=R(T)$
is the representation ring of $T$ and is  generated by formal exponents
$\{e^{\mu}\,:\,\mu\in P\}$, with the multiplication $e^{\lambda}\cdot
e^{\mu}=e^{\lambda+\mu}$.

Let $\omega$ be a dominant weight of $\mathfrak{g}$ and
$V(\omega)$ be an integrable irreducible highest weight representation
 of the  enveloping algebra $\CU(\mathfrak{g})$
with the highest weight $\omega$. For any   $w\in W$ the weight
$w(\omega)$ subspace $V^{[w(\omega)]}(\omega)$ in $V(\omega)$ is
one-dimensional. Let  $V_{w}(\omega)\subseteq V(\omega)$ be
a $\CU(\mathfrak{b})$-submodule generated by  enveloping algebra
$\CU(\mathfrak{b})$ of the Borel subalgebra $\mathfrak{b}$
 acting on $V^{[w(\omega)]}(\omega)$.
The $\CU(\mathfrak{b})$-module $V_{w}(\omega)$ is called Demazure
module. Characters of $V_{w}(\omega)$ are defined as
$${\rm ch}_{V_{w}(\omega)}=\sum_{\mu \in P}(dim
V^{[\mu]}_{w}(\omega))e^{\mu},
$$
and can be calculated using Demazure operators as follows. Define
Demazure operators corresponding to simple root $\alpha_i$ as
$$
\CD_{s_i}\,e^{\mu}=\frac{e^{\mu}-e^{-\alpha_i}e^{s_i(\mu)}}{1-e^{-\alpha_i}},
$$
where $s_i\in W$ is a simple reflection corresponding to $\alpha_i$.
Demazure operators  commute with $W$-invariant elements in
$\IZ[T]$  and satisfy the following relations
$$
\CD_{s_i}^2=\CD_{s_i}, \qquad (\CD_{s_i}\CD_{s_j})^{m_{ij}}=1,
$$
where $m_{ij}$ are equal to
$$
m_{ij}=2,\,3,\,\,4,\,\,6,\,\,\infty,
$$ for the values of entries of Cartan matrix $\|a_{ij}\|$ satisfying
$$
 a_{ij}a_{ji}=0,\,\,1,\,\,2,\,\,3,\,\,\geq 4,
$$
respectively. Here  we imply  $x^{\infty}:=1$. These relations
provide a correctly defined map $w\mapsto\CD_{w}$:
$$
w=s_{i_1}s_{i_2}\cdots s_{i_j}\longmapsto\CD_{w}=
\CD_{s_{i_1}}\CD_{s_{i_2}}\cdots \CD_{s_{i_{j}}}.
$$
Given a reduced (minimal length)  decomposition
$w=s_{i_1}s_{i_2}\cdots s_{i_j}$ of an element $w\in W$ we have for
the character of $V_{w}(\omega)$ \be {\rm
ch}_{V_{w}(\omega)}=\CD_{s_{i_1}}\CD_{s_{i_2}}\cdots
\CD_{s_{i_{j}}}\,e^{\omega}. \ee Now let $\mathfrak{g}$ be the
affine Lie algebra $\widehat{\mathfrak{gl}}_{\ell+1}$. The
corresponding root system can be realized as  a set of vectors in
$\IR^{\ell+2,1}$ supplied with a bi-linear symmetric form defined in
the bases $\{e_1,\ldots,e_{\ell+1},e_+,e_-\}$ by
$$
(e_i,e_j)=\delta_{ij},\qquad (e_{\pm},e_i)=(e_{\pm},e_{\pm})=0,\qquad
(e_{+},e_{-})=1.
$$
Simple roots of
$\widehat{\mathfrak{gl}}_{\ell+1}$ are given by
$$
\alpha_1=e_1-e_2,\qquad \alpha_2=e_2-e_3,\qquad \ldots \qquad
\alpha_{\ell}=e_{\ell}-e_{\ell+1}, \qquad
\alpha_0=e_+-(e_1-e_{\ell+1}).
$$
The fundamental weights $\omega_0,\omega_1\ldots,\omega_{\ell+1}$
are defined by the conditions
$(\omega_i,\alpha_j^{\vee})=\delta_{ij}$
$$
\omega_1=e_1+e_-,\qquad \omega_2=e_1+e_2+e_-,\qquad \ldots \qquad
\omega_{\ell+1}=\sum_{j=1}^{\ell+1} e_j+e_-,\qquad \omega_0=e_-.
$$
In the following we will also use the standard  notation
$\delta=\alpha_0+\sum_{i=1}^{\ell}\alpha_i=e_+$.
The Weyl group  $W$ has natural decomposition $W=\dot{W}\times
Q$ where $Q$ is a lattice generated by simple
coroots and  $\dot{W}$ is the  Weyl group  of the  finite-dimensional
Lie algebra $\mathfrak{gl}_{\ell+1}$.
Define a projection of the weight lattice $P$ of
$\widehat{\mathfrak{gl}}_{\ell+1}$ onto the  weight lattice
$\dot{P}$ of the finite-dimensional algebra $\mathfrak{gl}_{\ell+1}$
$$
\omega=\lambda_1e_1+\cdots +\lambda_{\ell+1}e_{\ell+1}+ke_-+re_+\longrightarrow
\dot{\omega}=\lambda_1e_1+\cdots +\lambda_{\ell+1}e_{\ell+1}.
$$
The projection on the lattice $\dot{P}$  of the action of the
generators of $W$  on $e_-+\sum \lambda_ie_i$ is given
by
$$s_i\,(\lambda_1,\ldots ,\lambda_{\ell+1})=(\lambda_1,\ldots
,\lambda_{i-1},
\lambda_{i+1},
\lambda_i,\lambda_{i+2},\ldots,\lambda_{\ell+1}),
$$
$$s_0\cdot(\lambda_1,\ldots
,\lambda_{\ell+1})=(\lambda_{\ell+1}+1,\lambda_2,\ldots ,
\lambda_{\ell},\lambda_{1}-1).
$$
Note that $|\omega|=\sum_{i=1}^{\ell+1}\,\lambda_i$ is invariant under the action
of $W$.
\begin{lem}
A set of orbits of $W$ acting on the weight lattice $\dot{P}$ of
$\mathfrak{gl}_{\ell+1}$ can be identified with $\IZ\times
(\IZ/\IZ_{\ell+1})$ and a set of representatives can be chosen as
follows
$$
\dot{\lambda}_{k,i}=(k+1,\ldots ,k+1,k,\ldots, k)=k\dot {\bf
1}+\dot{\omega_i},
$$
where ${\bf 1}=(1,\ldots,1)$ and $\dot{\omega}_i$ are fundamental
weights of $\mathfrak{gl}_{\ell+1}$.
\end{lem}

\noindent {\it Proof}:  Using $\dot{W}$
one can transform any weight of $\mathfrak{gl}_{\ell+1}$ to
a  dominant one $\dot{\lambda}=(\lambda_1\geq \lambda_2 \geq 
\cdots \geq \lambda_{\ell+1})$.
Now using elements $W$ transforming  dominant weights to dominant
we can change weights in such a way that the difference
$\lambda_{i}-\lambda_{i+1}$ is either $1$ or $0$ $\Box$

Define homomorphism
$$
\pi\,:\,\,\IZ[T] \rightarrow \IZ[z_1,\ldots, z_{\ell+1},q]
$$
$$
\pi(e^{\omega_i})=z_1\cdots z_i, \qquad \pi(e^{\omega_0})=1, \qquad
\pi(e^{\delta})=q.
$$
The following result was proved by Sanderson \cite{San1}.
\begin{te}
Let $\lambda_{k,i}=\omega_0+\dot{\lambda}_{k,i}$ and $\dot{\lambda}_{k,i}\in \dot{P}^+$
is given by
$$
\dot{\lambda}_{k,i}=(k+1,\ldots ,k+1,k,\ldots k)=k\cdot {\bf
1}+\dot{\omega}_i.
$$
Let $w\in W$ be such that for $\lambda=w\cdot \lambda_{k,i}$ its projection
$\dot{\lambda}$ be antidominant  weight i.e. $\lambda_1\leq \lambda_2 \leq \cdots
\leq \lambda_{\ell+1}$. Define $\dot{\lambda}'=w_0\dot{\lambda}$, where
$w_0\in\dot{W}$ is an element having a reduce decomposition
of maximal length.  Then the character of the Demazure module
$V_{w}(\lambda_{k,i})$ satisfy the following relation
$$
\pi\Big(ch_{V_{w}(\lambda_{k,i})} \Big)
=q^{\frac{1}{2}(\dot{\lambda},\dot{\lambda})-\frac{1}{2}(\dot{\lambda}_{k,i}, \dot{\lambda}_{k,i})}\,
P_{\dot{\lambda}'}(z;q,t=0),
$$
where $P_{\dot{\lambda}'}(z;q,t)$ is a Macdonald polynomial
corresponding to dominant partition $\dot{\lambda}'$ (see Definition
\ref{Macdef}).
\end{te}
The modified $q$-deformed
$\mathfrak{gl}_{\ell+1}$-Whittaker function is given  by 
$$\widetilde{\Psi}^{\mathfrak{gl}_{\ell+1}} _{z}(\underline{p}_{\ell+1})\,=\,
\Delta(\underline{p}_{\ell+1})\Psi^{\mathfrak{gl}_{\ell+1}}
_{z} (\underline{p}_{\ell+1}),$$
where
$$
\Delta(\underline{p}_{\ell+1})\,=\,\prod_{i=1}^{\ell}\,(p_{\ell+1,i}-
p_{\ell+1,i+1})_q!
$$
and $\Psi^{\mathfrak{gl}_{\ell+1}}_{z} (\underline{p}_{\ell+1})$
is defined by \eqref{main}.

\begin{te}\label{corSan} The following representation for the modified $q$-deformed
$\mathfrak{gl}_{\ell+1}$-Whittaker function holds
$$
\widetilde{\Psi}_{z}^{\mathfrak{gl}_{\ell+1}}(\underline{p}_{\ell+1})=
q^{\frac{1}{2}(\dot{\lambda}_{k,i},
\dot{\lambda}_{k,i})-\frac{1}{2}(\dot{\lambda},\dot{\lambda})}
\pi\Big(ch_{V_{w}(\lambda_{k,i})}\Big), \qquad
p_{\ell+1,i}=(\dot{\lambda}')_i.
$$
\end{te}
Thus the finite sum  \eqref{main} up to a simple multiplier
provides expression for a characters of affine Lie algebra Demazure module.

\noindent {\it Example.} Let us consider as an example the case of
$\ell=1$. We  have
$$
ch_{V_{(s_1s_0)^m}(\omega_0)}=\CD_{(s_1s_0)^m}\,e^{\omega_0},\qquad
ch_{V_{s_1(s_0s_1)^m}(\omega_1)}=\CD_{s_1(s_0s_1)^m}\,e^{\omega_1},
$$
where $\omega_0=e_-$ and $\omega_1=e_-+e_1$. Note that due
to the identity $\CD_1^2=\CD_1$ both characters are $W_1=S_2$
invariant and thus are given by linear combination of
$\mathfrak{gl}_2$-characters.
$$
 \dot{\lambda}_{0,0}=(0,0),\qquad \lambda=(s_1s_0)^m\omega_0, \qquad
 \dot{\lambda}=(-m,m),\qquad \dot{\lambda}'=(m,-m),
$$
$$
\dot{\lambda}_{0,1}=(1,0),\qquad \lambda=s_1(s_0s_1)^m\omega_1,\qquad
 \dot{\lambda}=(-m,m+1), \qquad \dot{\lambda}'=(m+1,-m),
$$
and thus
$$
\pi\Big(ch_{V_{(s_1s_0)^m(\omega_0)}}\Big)=
q^{m^2}\,P_{m,-m}(z_1,z_2;q,t=0),
$$
$$
\pi\Big(ch_{V_{s_1(s_0s_1)^m(\omega_1)}}\Big)
=q^{m(m+1)}\,P_{m+1,-m}(z_1,z_2;q,t=0).
$$

Let us note that there exists  a generalization of the results in \cite{San1} to the
case of simply-laced affine Lie algebras \cite{I}. Also the
structure of Demazure modules for arbitrary simply-laced affine
Lie group was clarified in  \cite{FL}.
It was shown that as a module
over corresponding finite-dimensional Lie algebra it is a finite
tensor product of finite-dimensional irreducible representations. In
the special case of $\widehat{\mathfrak{g}}=\widehat{\mathfrak{gl}}_{\ell+1}$
this is in complete agreement  with the  Proposition 3.4 in
\cite{GLO1}. Note that the case $\ell=1,2$ was considered before in
\cite{San2}.    All this seems implies that the connection between
$q$-deformed Whittaker functions, specialization of Macdonald
polynomials and Demazure modules discussed above  can be rather
straightforwardly generalized at least to the simply-laced case. We
conjecture that this indeed so and are going to discuss the details
elsewhere.

\section{$q$-deformed Whittaker function as a matrix element}

According to Kostant \cite{Ko} $\mathfrak{g}$-Whittaker functions can be understood
as matrix elements of infinite-dimensional representations of
universal enveloping algebra $\CU(\mathfrak{g})$ with action of Cartan
subalgebra $\mathfrak{h}\subset \mathfrak{g}$ integrated to the action
of the corresponding Cartan subgroup $H\subset G$. This interpretation
can be generalized to the case of $q$-deformed $\mathfrak{g}$-Whittaker functions
considered as matrix elements of infinite-dimensional representations
of quantum groups $\CU_q(\mathfrak{g})$ (see e.g. \cite{Et}).
 In this Section we derive a
representation of $q$-deformed $\mathfrak{gl}_{\ell+1}$-Whittaker
functions given explicitly by \eqref{main} as a matrix element of an
infinite-dimensional representation of multidimensional quantum
torus. Our construction is based on an iterative application of the
following standard identity  (see e.g. \cite{CK})
$$
(X+T)^n=\sum_{m=0}^n {n\choose m}_{\!\!q}\,X^{m}T^{n-m},\qquad TX=qXT
$$
where ${n\choose m}_{q}=(n)_q!/(m)_q!(n-m)_q!$ is a $q$-binomial coefficient.
This representation should arise in the Kostant framework
using a realization of $\CU_q(\mathfrak{gl}_{\ell+1})$ by difference
operators generalizing Gauss-Givental realization of
$\CU(\mathfrak{gl}_{\ell+1})$ proposed in \cite{GKLO}.  We check this
directly for $\mathfrak{g}=\mathfrak{sl}_2$ leaving the general case
to another occasion.

Let $\CA^{(\ell)}$ be an associative algebra,
$\{X^{\pm 1}_{k,i},T^{\pm 1}_{k,i}\}$, $k=1,\ldots,\ell$;
$i=1,\ldots,k$ be a complete set of generators  satisfying  relations
\be
 T_{k,i}X_{m,j}\,=\,q^{\delta_{k,m}\delta_{i,j}}\cdot X_{m,j}T_{k,i}.\ee
Introduce a set of polynomials
$f_{n,i}(z)\in\CA^{(\ell)}[z_1,\ldots,z_{\ell+1}]$
$$
f_{n,i}=f_{n,i}(\underline{z};X_{k,j},T_{k,j}),\qquad
n=1,\ldots,\ell+1;\quad i=1,\ldots,n,
$$
of degree ${\rm deg}\,f_{n,\,i}=i$ in variables
$\underline{z}=(z_1,\ldots,z_{\ell+1})$, defined by the following
recursive relations: \be f_{n,\,i}\,=\,f_{n-1,\,i}X_{n-1,\,i}\,+\,
z_n\,f_{n-1,\,i-1}T_{n-1,\,i},\hspace{1.5cm}i<n,\ee
where
$f_{n,0}=f_{00}=1$ and $f_{n,n}=z_1\cdots z_n$ with
$$
f_{n,n}\,=\,z_n\cdot f_{n-1,n-1}\hspace{2cm} n=1,\ldots,\ell+1.
$$
 In particular, we have $f_{11}=z_1$ and
$f_{21}=f_{11}X_{11}+z_2f_{10}T_{11}=z_1X_{11}+z_2T_{11}$.

Let $\CV$ be a representation of $\CA^{(\ell)}$, and $\CV^*$ be its
dual. Consider  $|v_+\>\in \CV$, $\<v_-|\in \CV^*$ such that
$\<v_-|v_+\>=1$ and  satisfying the conditions
$$
T_{i,k}|v_+\>=|v_+\>,\hspace{1.5cm}
\<v_-|X_{i,k}=\<v_-|,\hspace{1.5cm} i=1,\ldots,\ell,\,\,\,
k=1,\ldots ,i.
$$
Let us introduce normalized $q$-deformed
$\mathfrak{gl}_{\ell+1}$-Whittaker function as
$$
\widetilde{\Psi}^{\mathfrak{gl}_{\ell+1}}_{z_1,\ldots,z_{\ell+1}}
(\underline{p}_{\ell+1})\,=\,
\prod_{k=1}^\ell(p_{\ell+1,\,k}-p_{\ell+1,\,k+1})_q!\,\,
\Psi^{\mathfrak{gl}_{\ell+1}}_{z_1,\ldots,z_{\ell+1}}
(\underline{p}_{\ell+1}).
$$

\begin{te} The following representation of $q$-deformed
$\mathfrak{gl}_{\ell+1}$-Whittaker function holds
\be\label{NonCommutativeFormula}
\widetilde{\Psi}^{\mathfrak{gl}_{\ell+1}}_{z_1,\ldots,z_{\ell+1}}
(\underline{p}_{\ell+1})\,=\,\,\Big\<v_-\,\Big|\,
\prod_{k=1}^{\ell+1}\,f_{\ell+1,\,k}
(\underline{z};X_{n,i},T_{n,i})^{p_{\ell+1,\,k}-p_{\ell+1,\,k+1}}\,
\Big|\,v_+\,\Big\>, \ee where we assume $p_{\ell+1,\,\ell+2}=0$.
\end{te}

\noindent {\it Proof}:  Let us prove the Theorem by induction. We
assume that the representation \eqref{NonCommutativeFormula} for
$\mathfrak{gl}_{\ell}$ holds
$$
\widetilde{\Psi}^{\mathfrak{gl}_\ell}_{z_1,\ldots,z_\ell}
(\underline{p}_\ell)\,=\,\,\<v_-\,|\, \prod_{k=1}^\ell\,f_{\ell,\,k}
(\underline{z}')^{p_{\ell,\,k}-p_{\ell,\,k+1}}\,|\,v_+\,\>,
$$
where $\underline{z}'=(z_1,\ldots,z_\ell)$. The following recursive
relation follows from \eqref{qtodarec} \be\label{renrec}
\widetilde{\Psi}^{\mathfrak{gl}_{\ell+1}}_{z_1,\ldots,z_{\ell+1}}
(\underline{p}_{\ell+1})\,=\,
\sum_{\underline{p}_\ell\in\CP_{\ell+1,\ell}}\,z_{\ell+1}^{\sum
p_{\ell+1,i}-\sum p_{\ell,k}}\,\widetilde{Q}_{\ell+1,\ell}
(\underline{p}_{\ell+1},\underline{p}_\ell|q)\,
\widetilde{\Psi}^{\mathfrak{gl}_{\ell}}_{z_1,\ldots,z_{\ell}}
(\underline{p}_{\ell}), \ee
 where
$$
\widetilde{Q}_{\ell+1,\ell}
(\underline{p}_{\ell+1},\underline{p}_\ell|q)\,=\,
\prod_{k=1}^\ell{p_{\ell+1,k}-p_{\ell+1,k+1}\choose
p_{\ell,k}-p_{\ell+1,k+1}}_{\!\!q}
$$
Then \eqref{NonCommutativeFormula} for $\mathfrak{gl}_{\ell+1}$ is
obtained by repeated application of the identities
 \be
\sum_{p_{\ell,k}=p_{\ell+1,k+1}}^{p_{\ell+1,k}}\,\,
z_{\ell+1}^{p_{\ell+1,k}-p_{\ell,k}}
{p_{\ell+1,k}-p_{\ell+1,k+1}\choose
p_{\ell,k}-p_{\ell+1,k+1}}_{\!\!q}\,\cdot\\
\cdot\,\,(f_{\ell,k-1})^{p_{\ell,k-1}-p_{\ell,k}}
(f_{\ell,k})^{p_{\ell,k}-p_{\ell+1,k+1}} \,\,
X_{\ell,k}^{p_{\ell,k}-p_{\ell+1,k+1}}\,\,
T_{\ell,k}^{p_{\ell+1,k}-p_{\ell,k}}\,=\\
\nonumber =\,(f_{\ell,k-1})^{p_{\ell,k-1}-p_{\ell+1,k}}
\Big(f_{\ell,k}X_{\ell,k}\,+\,z_{\ell+1}f_{\ell,k-1}T_{\ell,k}
\Big)^{p_{\ell+1,k}-p_{\ell+1,k+1}} ,\ee to convert $q$-binomial
factors $\widetilde{Q}_{\ell+1,\ell}$ in \eqref{renrec} $\Box$

 The representation \eqref{NonCommutativeFormula} can be understood
as a particular realization of $q$-deformed
$\mathfrak{gl}_{\ell+1}$-Whittaker function as a matrix element of
an infinite-dimensional representation of
$\CU_q(\mathfrak{gl}_{\ell+1})$. In the following we demonstrate
this for the simplest case $q$-deformed $\mathfrak{sl}_2$-Whittaker
function.

Quantum deformed  universal enveloping algebra
$\CU_q(\mathfrak{sl}_2)$ is generated by generators $E,F,K$ satisfying the
relations
$$
KE=qEK,  \qquad KF=q^{-1}FK, \qquad EF- FE=-\frac{K-K^{-1}q}{1-q}.
$$
The center of $\CU_q(\mathfrak{sl}_2)$ is generated by a Casimir
element
$$
C=K+K^{-1}+(q+q^{-1}-2)FE.
$$
In irreducible representations
the image of $C$  is proportional to unite operator and we
parametrize the corresponding eigenvalue $c$  as follows
$$
c=-(z+z^{-1}).
$$
Consider a realization of $\CU_q(\mathfrak{sl}_2)$
 (see e.g. \cite{KLS})
$$
K=-zu^{-1},\qquad
E=\frac{v^{-1}(1-u^{-1})}{1-q},\qquad
F=-\frac{v(z-qz^{-1}u)}{1-q},
$$
where $uv=qvu$.  The general
$q$-deformed $\mathfrak{sl}_2$-Whittaker function
 is given by (compare with (2.17) in \cite{KLS}
 with $\omega_1=1$, $q=\exp(2\imath \pi  \omega_1/\omega_2)$)
\be\label{qGivsltwo} \Phi^{(\alpha_1,\alpha_2)}_z(x)=e^{-\pi x}\,q^{\imath x/2}
\,\, \<\psi_L^{(\alpha_1)}|\, q^{\imath \frac{x}{2}H}\,|\psi^{(\alpha_2)}_R\>,
\ee where $K=q^{H/2}$  and
$\psi^{(\alpha)}_{L}/\psi^{(\alpha)}_{R}$ are left/right Whittaker vectors
defined by
$$
E\psi^{(\alpha)}_L= q^{1-\alpha}\,e^{\imath \pi
\alpha}\frac{K^{\alpha}}
{1-q}\psi^{(\alpha)}_L,\qquad
F\psi^{(\alpha)}_R=e^{\imath \pi
\alpha}\frac{K^{-\alpha}}{1-q}\psi^{(\alpha)}_R,
$$
 \be\label{qTodaold}
(\CT^{-1}+\CT-q^{(\alpha_1-\alpha_2+1)}q^{-\imath x}
  \CT^{\alpha_1-\alpha_2})\Phi^{(\alpha_1,\alpha_2)}_{z}(x)\,=\,
(z+z^{-1}) \Phi^{(\alpha_1,\alpha_2)}_{z}(x), \ee where  $\CT\,f(x)
=f(x+\imath)$.

We would like to compare this representation with a
representation given in the previous section. The representation for
$\ell=1$ adopted to the case of $\mathfrak{sl}_2$
is given by
$$
\Psi_z(n)=\frac{1}{(n)_q!}\<v_-|\,(zX+z^{-1}T)^{n}\,|v_+\>,\qquad
\widetilde{\Psi}_z(n)=\<v_-|\,(zX+z^{-1}T)^{n}\,|v_+\>,
$$
where $T|v_+\>=|v_+>$ and $\<v_-|X=\<v_-|$. The functions
$\Psi_z(n)$ and $\widetilde{\Psi}_z(n)$ satisfy the following
equation \be\label{qTodanew}
(\CT^{-1}+(1-q^{n+1})\CT)\Psi_z(n)=(z+z^{-1})\Psi_z(n),\\
((1-q^n)\CT^{-1}+\CT)\widetilde{\Psi}_z(n)=(z+z^{-1})\widetilde{\Psi}_z(n),
\ee
where $\CT\,f(n) =f(n+1)$.
To reconcile the last equation in \eqref{qTodanew} and
the equation \eqref{qTodaold} we take
$$
\alpha_1=1,\,\,\,\,\alpha_2=2,\qquad x=\imath n.
$$
Then one has
$$
\Phi^{(1,2)}_z(\imath n)= \<\psi_L^{(1)}|\,(-K)^{-n}\,|\psi^{(2)}_R\>,
$$
and one  would like to have the following relation for the
$q$-deformed Whittaker function
\be\label{identification}
\widetilde{\Psi}_z(n)=\Phi^{(1,2)}_z(\imath n ).
\ee
To make this identification one should
should have the following relation
$$
-K^{-1}=z^{-1}u=z^{-1}T+zX.
$$
and demonstrate that $\<v_-|$ and $|v_+\>$ provide as representation
for left and right Whittaker vectors $\<\psi_L^{(1)}|$ and $|\psi^{(2)}_R\>$.

Consider the following unitary transformation
$$
U^{-1}(u,v)u U(u,v)=u+z^2v,\qquad U(u,v)=\prod_{n=0}^{\infty}(1+z^2vu^{-1}q^n)^{-1}\,=\,
\Gamma_q(-z^2vu^{-1}),
$$
where
$$
\Gamma_q(x)=\frac{1}{\prod_{j=0}^{\infty} (1-xq^j)}.
$$
Thus we have
$$
U^{-1}(u,v)\,v\,U(u,v)=(1+z^2vu^{-1})v,
$$
$$
U^{-1}(u,v)\,u\,U(u,v)=(1+z^2vu^{-1})u,
$$
$$
U^{-1}(u,v)\,v^{-1}\,U(u,v)=(1+q^{-1}z^2vu^{-1})^{-1}v^{-1},
$$
$$
U^{-1}(u,v)\,u^{-1}\,U(u,v)=(1+q^{-1}z^2vu^{-1})^{-1}u^{-1}.
$$
The conjugated generators are then given by
$$
\widehat{K}=U^{-1}(u,v)\,K\,U(u,v)=-z
(1+q^{-1}z^2vu^{-1})^{-1}u^{-1},
$$
$$
\widehat{E}=U^{-1}(u,v)\,E\,U(u,v)
=\frac{1}{1-q}\,{(1+q^{-1}z^2vu^{-1})^{-1}v^{-1}
(1-(1+q^{-1}z^2vu^{-1})^{-1}u^{-1})},
$$
$$
\widehat{F}=U^{-1}(u,v)\,F\,U(u,v)= -\frac{1}{1-q}
(1+z^2vu^{-1})v(z-z^{-1}q(1+z^2vu^{-1})u).
$$
\begin{prop} The following identities hold
$$
\widehat{E}|u=1\>=-\frac{1}{1-q}\,\widehat{K}^2|u=1\>,
\qquad
\widehat{F}|v=1\>=\frac{1}{1-q}\,\widehat{K}^{-1}|v=1\>,
$$
where $v|v=1>=|v=1>$, $u|u=1>=|u=1>$.
\end{prop}
{\it Proof}: Direct calculations $\Box$

Thus one can identify $|v_-\>\equiv |u=1\>=|\psi_L^{(2)}\>$,
$|v_+\>\equiv |v=1\>=|\psi^{(1)}_R\>$ in the
$U$-rotated bases. This also provides an identification \eqref{identification}.


\vskip 1cm

\noindent {\small {\bf A.G.}: {\sl Institute for Theoretical and
Experimental Physics, 117259, Moscow,  Russia; \hspace{8 cm}\,
\hphantom{xxx}  \hspace{4 mm} School of Mathematics, Trinity
College, Dublin 2, Ireland; \hspace{8 cm}\,
\hphantom{xxx}   \hspace{3 mm} Hamilton
Mathematics Institute, TCD, Dublin 2, Ireland;}}

\noindent{\small {\bf D.L.}: {\sl
 Institute for Theoretical and Experimental Physics,
117259, Moscow, Russia};\\
\hphantom{xxxxxx} {\it E-mail address}: {\tt lebedev@itep.ru}}\\

\noindent{\small {\bf S.O.}: {\sl
 Institute for Theoretical and Experimental Physics,
117259, Moscow, Russia};\\
\hphantom{xxxxxx} {\it E-mail address}: {\tt Sergey.Oblezin@itep.ru}}

\begin{thebibliography}{AB}
{
\frenchspacing \smallbreak

\bibitem[AOS]{AOS} H.~Awata, S.~Odake, J.~Shiraishi,
{\it Integral representations of the Macdonald symmetric functions},
Commun. Math. Phys. {\bf 179} (1996), 647--666; {\tt
arXiv:q-alg/9506006}.

\bibitem[CS]{CS}\, W.~Casselman, J.~Shalika, {\it
 The unramified principal series of p-adic groups II. The Whittaker
function}, Comp. Math. {\bf 41} (1980), 207--231.

\bibitem[Ch]{Ch} I.~Cherednik, {\it
Double Affine Hecke Algebras}, Cambridge University Press 2005.

\bibitem[CK]{CK}\, P.~Cheung, V.~Kac, {\it Quantum Calculus},
  Springer 2001.

\bibitem[De]{De} M.~Demazure, {\it
 D\'{e}singularisation des vari\'{e}t\'{e}s de Schubert
 g\'{e}n\'{e}ralis\'{e}es}.
 Ann. Sci. Ec. Norm. Sup\'{e}r. 7, 53-88 (1974).


\bibitem[Et]{Et}\, P.~Etingof, {\it Whittaker functions on quantum
groups and $q$-deformed Toda operators}, Amer. Math. Soc. Transl.
Ser.2, vol. 194, 9--25,  Amer.Math.Soc., Providence, RI,
1999;  {\tt arXiv:math.QA/9901053}.



\bibitem[FL]{FL} G.~Fourier, P.~Littelmann,  {\it Weyl
  modules, Demazure modules, KR-modules crystals, fusion products and
  limit constructions}, Adv. Math. {\bf 211} (2007), no. 2, 566-593;
{\tt arXiv:math.RT/0509276}.



\bibitem[GKLO]{GKLO} A.~Gerasimov, S.~Kharchev, D.~Lebedev, S.~Oblezin,
{\it On a Gauss-Givental representation of quantum Toda chain wave
function}, Int. Math. Res. Notices, (2006), AricleID 96489, 23
pages; {\tt arXiv:math.QA/0505310}.



\bibitem[GLO1]{GLO1} A.~Gerasimov,  D.~Lebedev,  S.~Oblezin, {\it
 On q-deformed $\mathfrak{gl}_{\ell+1}$-Whittaker functions  I},
 {\tt arXiv:math.RT/0803.0145}.


\bibitem[GLO2]{GLO2} A.~Gerasimov,  D.~Lebedev,  S.~Oblezin, {\it
 On q-deformed $\mathfrak{gl}_{\ell+1}$-Whittaker functions  II},
 {\tt arXiv:math.RT/0803.0970}.

\bibitem[GLO3]{GLO3} A.~Gerasimov,  D.~Lebedev,  S.~Oblezin, {\it
 On q-deformed $\mathfrak{gl}_{\ell+1}$-Whittaker functions  IV},
 to appear.


\bibitem[I]{I} B.~Ion, {\it Nonsymmetric Macdonald polynomials and
    Demazure characters}, Duke Mathematical Journal 116 (2003),  no. 2,  299-318.

\bibitem[KL]{KL}\,
 S.~Kharchev, D.~Lebedev,
{\it Eigenfunctions of $GL(N,R)$ Toda chain:  The Mellin-Barnes
representation}, JETP Lett. {\bf 71} (2000), 235--238; {\tt
arXiv:hep-th/0004065}.


\bibitem[KLS]{KLS}\, S.~Kharchev, D.~Lebedev,
M.~Semenov-Tian-Shansky, {\it Unitary representations of
$U_q(sl(2,R))$, the modular double and the multiparticle q-deformed
Toda chains}, Comm. Math. Phys. {\bf 225} (2002), 573--609; {\tt
arXiv:hep-th/0102180}.

\bibitem[Ko]{Ko}
B. Kostant, {\it Quantization and representation
theory}, In: Representation Theory of Lie Groups, Proc. of Symp.,
Oxford, 1977, pp. 287-317, London Math. Soc. Lecture Notes
series, {\bf 34}, Cambridge, 1979.


\bibitem[Ku]{Ku}  S.~Kumar, {\it
Demazure character formula in arbitrary Kac-Moody setting},
 Invent. Math. {\bf 89} (1987), 395--423.

\bibitem[Mac]{Mac} I.G.~Macdonald, {\it
A new class of symmetric functions},
 S\'{e}minaire Lotharingien de Combinatoire, B20a (1988), 41 pp.

\bibitem[M]{M} O.~Mathieu, {\it Formules de charact\`{e}res pour les
    alg\`{e}bres Kac-Moody g\`{e}n\`{e}rales}, Ast\`{e}risque \\{\bf
    159-160}, Soc. Math. France, Monrouge, 1988.

\bibitem[San1]{San1} Y.~Sanderson, {\it On the Connection Between
    Macdonald Polynomials and Demazure Characters}, J. of Algebraic
  Combinatorics, {\bf 11} (2000), 269-275.

\bibitem[San2]{San2} Y.~Sanderson, {\it Real characters for Demazure
 modules of rank two affine Lie algebras}, J. Algebra {\bf 184} (1996).

\bibitem[Sh]{Sh}\, T.~Shintani, {\it On an
explicit formula for class 1 Whittaker functions on $GL_n$ over
p-adic fields}, Proc. Japan Acad. {\bf 52} (1976), 180--182.

}
\end{thebibliography}
\end{document}